\documentclass[12pt]{article}
\usepackage{amsmath}
\usepackage{amsfonts}
\usepackage{amssymb}
\usepackage{theorem}

\title{An $\omega$-categorical structure with amenable automorphism group}
\author{A. Ivanov
\thanks{The author is supported by Polish National Science Centre grant DEC2011/01/B/ST1/01406} } 
\date{ } 

\newtheorem{thm}{Theorem}[section] 
\newtheorem{lem}[thm]{Lemma}
\newtheorem{definicja}[thm]{Definition}
\newtheorem{cor}[thm]{Corollary}
\newtheorem{prop}[thm]{Proposition} 
\newtheorem{remark}[thm]{Remark} 

\newtheorem{example}{Example}[section] 
\begin{document}
\maketitle
\topskip 20pt

\begin{quote}
{\bf Abstract.} 
%\begin{abstract} 
We analyse $\omega$-categorical precompact 
expansions of particular $\omega$-categorical 
structures from the viewpoint of amenability of 
their automorphism groups. 
%\footnote{The author is supported by Polish National Science Centre grant DEC2011/01/B/ST1/01406} 
%\end{abstract}

\bigskip

{\em 2010 Mathematics Subject Classification:}  03C15, 03E15 

{\em Keywords:}  Amenable groups, Countably categorical structures.
\end{quote}

\bigskip

%%%%%%%%%%%%%%%%%%%%

\section{Introduction} 

A group $G$ is called {\bf amenable} if every $G$-flow 
(i.e. a compact Hausdorff space along with a continuous G-action) 
supports an invariant Borel probability measure.  
If  every $G$-flow has a fixed point then we say that $G$ 
is {\bf extremely amenable}. 
Let $M$ be a relational countably categorical structure which is 
a Fra\"{i}ss\'{e} limit of a Fra\"{i}ss\'{e} class $\mathcal{K}$. 
In particular $\mathcal{K}$ coincides with $Age(M)$, 
the class of all finite substructures of $M$. 
By Theorem 4.8 of \cite{KPT} the group $Aut(M)$ is 
{\em extremely amenable if and only if the class $\mathcal{K}$ 
has the Ramsey property and consists of rigid elements.} 
Here the class $\mathcal{K}$ is said to have 
the Ramsey property if  for any $k$ 
and a pair $A<B$ from $\mathcal{K}$ 
there exists $C\in \mathcal{K}$ so that 
each $k$-coloring 
$$
\xi :{C\choose A}\rightarrow k
$$ 
is monochromatic on some ${B'\choose A'}$
from $C$ which is a copy of 
${B\choose A}$, i.e. 
$$ 
C\rightarrow (B)^{A}_k . 
$$ 
We remind the reader that ${C \choose A}$ 
denotes the set of all substructures of $C$ 
isomorphic to $A$. 
This result has become a basic tool  to 
amenability of automorphism groups.  
To see whether $Aut(M)$ is amenable 
one usually looks for  an expansion $M^*$ 
of $M$ so that $M^*$ is a Fra\"{i}ss\'{e} structure 
with extremely amenable $Aut(M^* )$. 
Moreover it is usually assumed that $M^*$ is 
a {\bf precompact} expansion of $M$, 
i.e. every member of $\mathcal{K}$ has 
finitely many expansions in $Age(M^* )$, 
see \cite{KPT}, \cite{KS}, \cite{The}, \cite{AKL} 
and \cite{Z}.  
Theorem 9.2 from \cite{AKL} and Theorem 2.1 
from \cite{Z} describe amenability of $Aut(M)$ in 
this situation. 
The question if there is a countably 
categorical structure $M$ with amenable 
automorphism group which does not have 
expansions as above was formulated by 
several people. 
We mention very similar Problems 27, 28 
in \cite{BPT} where precompactness 
is replaced by $\omega$-categoricity 
and finite homogenity. 
\parskip0pt 

We think that in order to construct a required example 
one can use the ideas applied in \cite{ivanov} 
where we construct an $\omega$-categorical structure 
so that its theory is not G-compact 
and it does not have AZ-enumerations. 
These ideas develop ones applied in slightly different forms  
in \cite{ivma} and \cite{ivanov2} for some other questions. 
Moreover Casanovas, Pelaez and Ziegler suggest 
in \cite{cpz} a general method which simplifies and generalises 
our approach from \cite{ivanov2}, \cite{ivanov} and \cite{ivma}. 
The basic object of this construction is 
a  particular theory $T_E$ of equivalence 
relations $E_n$ on $n$-tuples.   
The paper \cite{cpz} pays attention 
to several model-theoretic properties 
of $T_E$. 

Below we study $T_E$ from the viewpoint of 
(extreme) amenability of its expansions. 
Then we apply our results to a construction of 
a family of concrete candidates for 
an example of an $\omega$-categorical 
structure with amenable automorphism group and 
without $\omega$-categorical precompact expansions with 
extremely amenable automorphism groups.  
We will in particular show that these structures have 
the following unusual combination of properties: 
\begin{itemize} 
\item the automorphis group is amenable; 
\item it does not satisfy Hrushovski's extension property;  
\item it does not have an order expansion with the Ramsey property.
\end{itemize}  
In fact we will show a slightly stronger version of the latter 
property.

\section{Equivalence relations} 

We start with a very interesting reduct of the structure from \cite{ivanov}.  
This is $T_E$ mentioned in the introduction. 
It has already deserved some attention in model-theoretic community, 
see \cite{cpz}.  

Let $L_0 =\{ E_n : 0< n<\omega \}$ be a first-order language,
where each $E_n$ is a relational symbol of arity $2n$.
Let ${\mathcal{K}}_0$ be the class of all finite $L_0$-structures $C$ 
where each relation $E_n (\bar{x},\bar{y})$ determines an equivalence
relation on the set (denoted by ${C \choose n}$) of unordered
$n$-element subsets of $C$. 
In particular for every $n$ the class ${\mathcal{K}}_0$ satisfies the sentence 
$$
\forall \bar{x}\bar{y} (E_n (x_1 ,...,x_n ,y_1 ,...,y_n )\rightarrow
\bigwedge \{ E_{n}(y_1 ,...,y_n ,x_{\sigma (1)},...,x_{\sigma (n)}):
\sigma \in Sym (n)\}).
$$ 
Note that  for $C\in {\mathcal{K}}_0$,  
$E_n$ is not satisfied by $\bar{a} ,\bar{b}$ if one of these 
tuples has a repetition. 
Thus for $n>|C|$ we put that no $2n$-tuple from $C$
satisfies $E_n (\bar{x},\bar{y})$.  
It is easy to see that ${\mathcal{K}}_0$ is closed under taking
substructures and the number of isomorphism types of
${\mathcal{K}}_0$-structures of any size is finite.
\parskip0pt

Let us verify {\em the amalgamation property for} ${\mathcal{K}}_0$.  
Given $A,B_{1},B_{2} \in {\mathcal{K}}_0$ with 
$B_{1}\cap B_{2}=A$, define $C\in {\mathcal{K}}_0$ 
as $B_{1}\cup B_{2}$ with the finest equivalence relations  
among those which obey the following rules.
When $n\le |B_{1}\cup B_{2}|$ and
$\bar{a}\in {B_1 \choose n} \cup {B_2 \choose n}$ we put
that the $E_n$-class of $\bar{a}$ in $C$ is contained in
${B_1 \choose n} \cup {B_2 \choose n}$.
We also assume that all $n$-tuples meeting both
$B_{1}\setminus B_{2}$ and $B_{2}\setminus B_{1}$
are pairwise equivalent with respect to $E_{n}$.
In particular if $n\ge max(|B_{1}|,|B_{2}|)$ we put that all
$n$-element $n$-tuples from $C$ are pairwise $E_{n}$-equivalent.
\parskip0pt

It is easy to see that this amalgamation also works 
for the joint embedding property.\parskip0pt

Let $M_0$ be the countable universal homogeneous structure 
for ${\mathcal{K}}_0$. 
It is clear that in $M_0$ each $E_n$ defines infinitely 
many classes and each $E_n$-class is infinite. 
Let $T_E = Th(M_0 )$.

Theorem \ref{nR} which we prove below, shows 
that $M_0$ cannot be treated by the methods of 
\cite{KPT}. 
It states that the group $Aut(M_0 )$ is amenable but 
the structure $M_0$ does not have a linear ordering 
so that the corresponding age has the order property 
and the Ramsey property.  

It is worth noting that this statement already holds for 
the $\{ E_1 ,E_2\}$-reduct of $M_0$, see the proof below. 
Thus our theorem also gives some interesting 
finitely homogeneous examples. 
On the other hand amenability of 
$Aut(M_0 )$ is a harder task than 
the corresponding statement in 
the reduct's case. 

The statement that $Aut(M_0 )$ is amenable 
is a consequence of a stronger property, 
namely {\em Hrushovski's extension property} for 
partial isomorphisms. 
This is defined for Fra\"{i}ss\'{e} limits 
as follows. 

\begin{definicja} 
A universal ultrahomogeneous structure $U$ 
satisfies Hrushovski's extension property 
if for any finite family of finite partial 
isomorphisms between substructures of $U$ 
there is a finite substructure 
$F< U$ containing these substructures so 
that any isomorphism from 
the family extends to an automorphism of $F$. 
\end{definicja} 

Proposition 6.4 of \cite{kechros} states that 
the structure $U$ has Hrushovski's 
extension property if and only if  
$Aut(U)$ has a dense subgroup which is the union 
of a countable chain of compact subgroups. 
The latter implies amenability by Theorem 449C
of \cite{fremlin}.

\begin{thm} \label{nR} 
(a) The structure $M_0$ satisfies Hrushovski's 
extension property. 
In particular the group $Aut(M_0 )$ is amenable. 

(b) 
The structure $M_0$ does not have any expansion by a linear order 
so that $Th(M_0 ,<)$ admits elimination of quantifiers 
and the age of $(M_0, <)$ satisfies the Ramsey property. 
\end{thm} 

The proof uses some material from \cite{herwig}. 
We now describe it.  

Let $L$ be a finite relational language. 
We say that an $L$-structure $F$ is {\em irreflexive} 
if for any $R\in L$, any tuple from $F$ satisfying $R$ 
consists of pairwise distinct elements. 
An irreflexive $L$-structure $F$ is called 
a {\em link structure} if $F$ is a singleton 
or $F$ can be enumerated  $\{ a_1 ,...,a_n \}$ 
so that $( a_1 ,...,a_n )$ 
satisfies a relation from $L$. 

Let ${\cal S}$ be a finite set of link structures. 
Then an $L$-structure $N$ is of {\em link type} 
${\cal S}$ if any substructure of $N$ which 
is a link structure is isomorphic 
to a structure from ${\cal S}$.  

An $L$-structure $F$ is {\em packed} if 
any pair from $F$ belongs to a link 
structure which is a substructure of $F$. 

If ${\cal R}$ is a finite family of packed 
irreflexive $L$-structures, then an 
$L$-structure $F$ is called ${\cal R}$-{\em free} 
if there does not exist a weak homomorphism 
(a map preserving the predicates) 
from a structure from ${\cal R}$ to $F$.  

Proposition 4 and Theorem 5 of \cite{herwig} state 
that for any family of irreflexive link structures 
${\cal S}$ and any finite family of irreflexive 
packed $L$-structures ${\cal R}$ the class of 
all irreflexive finite $L$-structures of link 
type ${\cal S}$ which are ${\cal R}$-free, 
has the free amalgamation property and 
Hrushovski's extension property  
for partial isomorphisms. 

We will use a slightly stronger version of 
this statement concerning {\em permorphisms}.  
A partial mapping $\rho$ on $U$ is called 
a $\chi$-{\bf permorphism}, if $\chi$ is 
a permutation of symbols in $L$ preserving 
the arity and for every $R\in L$ and 
$\bar{a}\in Dom(\rho )$ we have 
$$ 
\bar{a}\in R \Leftrightarrow \rho (\bar{a}) \in R^{\chi}.  
$$ 
The following statement is a version of 
Lemma 6 from \cite{herwig}. 

\begin{lem} \label{permorphism} 
Let $L$ be a finite language,  
$\chi_1 ,...,\chi_n$ be arity preserving 
permutations of $L$ and $\mathcal{S}$ be 
a finite $\{ \chi_i \}_{i\le n}$-invariant 
family of irreflexive link structures. 
Let $\mathcal{R}$ be a finite family of 
finite irreflexive packed $L$-structures 
of link type $\mathcal{S}$ so that 
$\mathcal{R}$ is invariant under all $\chi_i$. 
Let $A$ be a finite structure which belongs to 
the class, say $K$, of $L$-structures of 
link type $\mathcal{S}$ which are $\mathcal{R}$-free. 
Let $\rho_i$, $i\le n$, be parial $\chi_i$-permorphisms 
of $A$. 

Then there is a finite $B\in K$ containing $A$ 
so that each $\rho_i$ extends to a permutation of $B$ 
which is a $\chi_i$-permorphism.  
\end{lem} 

\bigskip

{\em Proof of Theorem \ref{nR}.} 
(a) 
For each $n>0$ enumerate all $E_n$-classes. 
Consider the expansion of $M_0$ by distinguishing 
each $E_n$-class by a predicate $P_{n,i}$ according 
the enumeration. 
Let $L^*$ be the language of all predicates $P_{n,i}$ 
and let $M^*$ be the $L^*$-structure defined on $M_0$. 
For every finite sublanguage $L'\subseteq L^*$  
let $M^* (L')$ be                                        
the $L'$-reduct of $M^*$ defined by these interpretations.  

We denote by $\mathcal{K}(L')$ the class of all finite 
$L'$-structures with the properties that for any arity $l$ 
represented by $L'$: 
\begin{itemize} 
\item any $l$-relation is irreflexive and invariant with respect to all permutations of variables,  
\item any two relations of $L'$ of arity $l$ have empty intersection.  
\end{itemize} 
Let ${\cal S}(L')$ be the set of all link structures of ${\cal K}(L')$ 
satisfying these two properties. 
Thus ${\cal K}(L')$ is of link type ${\cal S}(L')$. 

{\bf Claim 1.} For every finite sublanguage $L'\subseteq L^*$ 
the structure $M^* (L')$ is a universal structure 
with respect to the class $\mathcal{K}(L')$. 

It is easy to see that any structure from $\mathcal{K}(L')$ 
can be expanded to a structure from $\mathcal{K}_0$ 
so that $L'$-predicates become 
classes of appropriate $E_n$'s. 

{\bf Claim 2.} For every finite sublanguage $L'\subseteq L^*$ 
the structure $M^* (L')$ is an ultrahomogeneous structure. 

Let $f$ be an isomorphism between 
finite substructures of $M^* (L')$.  
We may assume that $Dom (f)$ contains 
tuples representing all $M^* (L')$-predicates 
of $L'$ (some disjoint tuples can be added 
to $Dom(f)$ in a suitable way).  
Then $f$ extends to an automorphism 
of $M_0$ fixing the classes of appropriate 
$E_n$'s which appear in $L'$. 
Thus this automorphism is 
an automorphism of $M^* (L')$ too. 

{\bf Claim 3.} 
For each finite sublanguage $L'\subseteq L^*$ 
let ${\cal R}(L')$ be the family of all packed $L'$-structures 
of the form $(\{ a_1 ,...,a_n \}, P_{n,i}, P_{n, j})$, 
where $i \not=j$, $P_{n, i}= \{(a_1, ...,a_n )\}$ and  
$P_{n,j}= \{(a_{\sigma(1 )}, ...,a_{\sigma (n )})\}$ for 
some permutation $\sigma$. 
Then the class $\mathcal{K}(L')$ is the class of all irreflexive 
finite $L'$-structures of link type ${\cal S}(L')$, 
which are ${\cal R}(L')$-free. 

The claim is obvious. 
By Proposition 4 and Theorem 5 of \cite{herwig} we now see that 
$\mathcal{K}(L')$ is closed under substructures, 
has the joint embedding property, the free amalgamation property,  
Hrushovski's extension property and its version for permorphisms, 
i.e. the statement of Lemma \ref{permorphism}. 

By Claim 1 and Claim 2 the structure $M^* (L')$ is 
the universal homogeneous structure of ${\cal K}(L')$. 
In particular any tuple of finite partial isomorphisms 
(permorphisms) of $M^* (L')$ can be extended to 
a tuple of automorphisms (permorphisms) of 
a finite substructure of $M^* (L')$. 
\parskip0pt 

Note that the same statement holds for the structure $M^*$. 
To see this take any tuple $f_1 ,...,f_k$ of finite partial 
isomorphisms (resp. $\chi_i$-permorphisms) of $M^*$ 
(assuming that $\chi_i$ are finitary). 
Let $r$ be the size of the union $\bigcup_{i\le k} Dom(f_i )$ 
and $L'$ be the minimal (resp. $\{ \chi_i \}_{i\le k}$-invariant) 
sublanguage of $L^*$ of arity $r$ containing of all 
relations of $M^*$ which meet any tuple from $\bigcup_{i\le k} Dom(f_i )$. 
Then there is a finite substructure $A$ of $M^* (L')$ containing 
$\bigcup_{i\le k} Dom(f_i )$ so that each $f_i$ extends to an automorphism 
(resp. $\chi_i$-permorphism) of $A$. 

Let $r'$ be the size of $A$.
Let $L''$ be a sublanguage of $L^*$ so that $L'\subseteq L''$
and for each arity $l\le r'$ the sublanguage $L''\setminus L'$ 
contains exactly one $l$-relation, say $P_{l,n_l}$ 
(fixed by $\{ \chi_i \}_{i\le k}$). 
Since $M^*$ is the universal homogeneous structure 
of ${\cal K}(L'')$ the substructure $A$ can be chosen so 
that any $l$-subset of $A$ which does not satisfy 
any relation from $L'$, does satisfy $P_{l, n_l}$. 

As a result any automorphism (permorphism) of $A$ 
preserves the relations of $L'''\setminus L'$  
for any $L^{'''}\subset L^*$ containing $L''$. 
Thus it extends to an automorphism (permorphism) 
of $M^* (L^{'''})$. 
In paricular it extends to an automorphism (permorphism) of $M^*$. 

As in Proposition 6.4 of \cite{kechros} we see that 
$Aut(M^* )$ has a dense subgroup which is the union 
of a countable chain of compact subgroups. 
In particlar we arrive at the following statement. 

{\bf Claim 4.} $Aut(M^* )$ is amenable. 

Since each automorphism of $M_0$ 
is a permorphism of $M^*$ and vice versa, 
we also see that 
$Aut(M_0  )$ has a dense subgroup which is the union 
of a countable chain of compact subgroups. 
In particular $Aut(M_0 )$ is amenable. 
%Since each automorphism of $M_0$ preserves all $E_i$, 
%$i>0$, it is easy to see that there is a natural 
%homomorphism from $Aut(M_0 )$ to the product 
%$\prod_{i>0} Sym ({M_0 \choose i }/E_i )$  
%and $Aut(M^* )$ is the kernel of it. 
%By Lemma \ref{support} this homomorphism 
%is surjective. 
%By Theorem 449C of \cite{fremlin} 
%we conclude that the group $Aut(M_0 )$ is amenable. 

(b) 
Consider a linearly ordered expansion $(M_0 ,<)$ together 
with the corresponding age, say $\mathcal{K}^<$. 
Assume that $\mathcal{K}^<$ has the Ramsey property. 

Note that  $\mathcal{K}^<$ does not contain 
any three-element structure of the form 
$a<b<c$, where $a$ and $c$ belong to the same 
$E_1$-class which is distinct from the $E_1$-class 
of $b$.  
Indeed, otherwise repeating the argument of Theorem 6.4 
from \cite{KPT}, we see that in any larger structure 
from $\mathcal{K}^<$ we can colour  two-elements structures 
$a<b$ with $\neg E_1 (a,b)$, so that there is no monochromatic 
three-element structure of the form above. 

As a result we see that any $E_1$-class of $(M_0 ,<)$ is convex. 
We now claim that the following structure $B$ 
can be embedded into $(M_0 ,<)$. 
  
Let $B=\{ a_1 <a_2 <a_3 <a_4 < b_1 <b_2 \}$, where 
the $E_1$-classes of all elements are pairwise distinct, 
but the pairs $\{ a_1 ,a_2 \}$ and $\{ b_1 ,b_2 \}$ 
are $E_2$-equivalent. 
We assume that in all other cases 
any two distinct pairs from $B$ belong 
to distinct $E_2$-classes. 
Moreover we assume that for each 
$k=3,4,5$ all $k$-subsets 
from $B$ belong to the same $E_k$-class. 
In particular the ordered structures 
defined on  $\{  a_1 ,a_2 ,a_3 ,a_4 \}$ 
and  $\{ a_3 ,a_4 ,b_1 ,b_2 \}$ 
are isomorphic. 
Let $A$ represent  
this isomorphism class.  

Since $M_0$ is the universal homogeneous structure 
with respect to $\mathcal{K}_0$, 
taking any tuple $a'_1 <a'_2 <a'_3 <a'_4 < b'_1 <b'_2$  
with pairwise distinct $E_1$-classes we can find  $B$  
in $M_0$ as a half of a copy of a structure from 
$\mathcal{K}_0$ consisting of 12 elements 
where each $E_1$-class is represented by a pair 
$(a'_i ,a_i )$ or $(b'_i ,b_i )$.  

To show that the Ramsey property does not 
hold for the age of $(M_0 ,<)$ take any finite 
substructure $C$ of this age which extends $B$. 
Fix any enumeration of $E_2$-classes 
ocurring in $C$. 
Then colour a copy of $A$ red if 
the class of the first two elements 
is enumerated before the class 
of the last pair. 
Otherwise colour such a copy green. 
It is clear that $C$ does not contain 
a structure isomorphic to $B$ so that all 
substructures of  type $A$ are of the same colour. 
$\Box$ 

\bigskip

\begin{remark} 
{\em It is worth noting that 
the class $\mathcal{K}^{<}_0$ of all linearly 
ordered members of $\mathcal{K}_0$ 
has JEP and AP,  i.e. there is 
a generic expansion of $M_0$ 
by a linear ordering. 
To see AP we just apply the amalgamation 
described above together with the 
standard amalgamation of orderings. }
\end{remark}

\section{Adding dense linear orders}

In order to obtain a structure with 
the properties as in Section 1, but without 
Hrushovski's extension property 
we use a general approach from \cite{cpz}. 
In fact our starting point is Corollary 2.8 from 
\cite{cpz} that sets ${M_0 \choose n}/E_n$ 
(definable in $Th^{eq} (M_0 )$) 
are stably embedded in $M_0$.  

We remind the reader that a $0$-definable 
predicate $P$ of a theory $T$ is called {\bf stably embedded} 
if every definable relation on $P$ is definable 
with parameters from $P$.   
If $M$ is a saturated model of $T$ then 
$P$ is stably embedded if and only if every 
elementary permutation of $P(M)$ extends to 
an automorphism of $M$ 
(see remarks after Definition 2.4 in \cite{cpz}). 
We now formulate Lemma 3.1 from \cite{cpz}. 
\begin{quote} 
Let $T$ be a complete theory with two sorts $S_0$ and $S_1$. 
Let $\tilde{T}_1$ be a complete expansion of 
$T\upharpoonright S_1$. 
Assume that $S_1$ is stably embedded. 
Then\\ 
(1) $\tilde{T} = T\cup \tilde{T}_1$ is a complete theory; \\ 
(2) $S_1$ is stably embedded in $\tilde{T}$ and 
$\tilde{T}\upharpoonright S_1 = \tilde{T}_1$. \\ 
(3) if $T$ and $\tilde{T}_1$ are $\omega$-categorical, 
then $\tilde{T}$ is also $\omega$-categorical. 
\end{quote} 

We now describe our {\bf variations} of $M_0$. 
Let us fix $S_n = {M_0 \choose n}/E_n$, 
$n\in \omega$, and consider them as a sequence 
of stably embedded sorts in $Th^{eq}(M_0 )$ 
(this is Corollary 2.8 of \cite{cpz}). 
We can distinguish relations 
$\{ a_1 ,..,a_n \} \in e$, where $e\in S_n$ 
is an $E_n$-class, $n\in \omega$.  

We also fix a subset $P\subset \omega \setminus \{ 1, 2 \}$  
and consider the language  
$$
L^S_P =\{ E_{n}: 0<n\in\omega \} \cup \{ S_n ,<_{S_n} : n\in P\} , 
$$   
where $<_{S_n}$ are binary relations on $S_n$.  
Let $\tilde{T}_1$ be the theory of 
sorts $\{ S_n : n\in\omega \}$, where for 
every $n\in P$ the relation $<_{S_n}$ is a dense linear 
order without ends. 
When $n\not\in P$ the sort $S_n$ is considered as 
a pure set. 
This is an $\omega$-categorical theory for each $S_n$. 
Applying Lemma 3.1 from \cite{cpz} 
we define the complete theory 
$T^S_P = T_E \cup \tilde{T}_1$
which is $\omega$-categorical and 
every sort $S_n$ is stably embedded into 
$T^S_P$.

We now define an one-sorted version of $T^S_P$. 
Its countable model will be the example 
anounced in Introduction. 
 
Let $L_P =\{ E_{n}: 0<n\in\omega \} \cup \{ <_n : n\in P\}$ 
be a first-order language, where each $E_{n}$ and 
$<_n$ is a relational symbol of arity $2n$.
The $L_P$-structure $M$ is built by 
the Fra\"{i}ss\'{e}'s construction. 
Let us specify a class ${\mathcal{K}}_P$ of 
finite $L_P$-structures, which will become 
the class of all finite substructures of $M$.
\parskip0pt 

Assume that in each $C\in {\mathcal{K}}_P$ each 
relation $E_{n}(\bar{x},\bar{y})$ determines 
an equivalence relation on the set (denoted by
${C \choose n}$) of unordered $n$-element subsets of $C$. 
As before for $C\in {\mathcal{K}}_P$ and $n>|C|$ we put that 
no $2n$-tuple from $C$ satisfies $E_n (\bar{x},\bar{y})$. 

For $n\in P$ the relations $<_n$ are irreflexive and respect $E_{n}$,
$$
\forall \bar{x},\bar{y},\bar{u},\bar{w} (E_{n}(\bar{x},\bar{y})
\wedge E_{n}(\bar{u},\bar{w})\wedge <_n(\bar{x},\bar{u})
\rightarrow <_n(\bar{y},\bar{w})). 
$$ 
Every $<_n$ is interpreted by a linear order 
on the set of $E_{n}$-classes.
Therefore we take the corresponding axioms 
(assuming below that tuples consist of pairwise distinct elements):    
$$
\forall \bar{x},\bar{y} (<_n(\bar{x},\bar{y})\rightarrow \neg E_{n}(\bar{x},\bar{y}) ); 
$$ 
$$
\forall \bar{x},\bar{y},\bar{z} (<_n(\bar{x},\bar{y}) \wedge 
<_n(\bar{y},\bar{z}) \rightarrow <_n(\bar{x},\bar{z}));
$$
$$
\forall \bar{x},\bar{y} (\neg E_{n}(\bar{x},\bar{y}) 
\rightarrow <_n(\bar{x},\bar{y}) \vee <_n (\bar{y},\bar{x})).
$$

\begin{lem} \label{LEM} 
(1) The class ${\mathcal{K}}_P$ satisfies the joint 
embedding property and the amalgamation property. \\ 
(2)  Let $M$ be the generic structure of 
${\mathcal{K}}_P$.  
For every $n>0$ let $M_n ={M\choose n }/E_n$.  \\ 
Then $Th(M)$ is $\omega$-categorical, 
admits elimination of quantifiers, and $<_n$ is a dense 
linear ordering on $M_n$ without ends (when $n\in P$). 
The structure $M$ is an expansion of $M_0$. \\ 
(3) Let $\rho_i$, $i\le k$, 
be a sequence of finitary maps on $M_i $ 
which respect $<_i$ for $i\in P$. 
Then there is an automorphism $\alpha \in Aut(M)$ 
realising each $\rho_i$ on its domain.   
\end{lem} 

{\em Proof.} 
(1) Given $A,B_{1},B_{2} \in \mathcal{C}$ 
with $B_{1}\cap B_{2}=A$, define 
$C\in \mathcal{K}$ as $B_{1}\cup B_{2}$.
The relations $E_{n},<_n$,
$n\le |B_{1}\cup B_{2}|$, are defined so 
that $C\in {\mathcal{K}}$, $B_1 <C$, $B_2 <C$ 
and the following conditions hold.
Let $n\le |B_{1}\cup B_{2}|$.
We put that all $n$-element $n$-tuples meeting both 
$B_{1}\setminus B_{2}$ and $B_{2}\setminus B_{1}$ 
are pairwise equivalent with respect to $E_{n}$.
We additionally demand that they are 
equivalent to some tuple from some 
$B_{i}$, $i\in \{ 1,2\}$, 
if $n\le max(|B_{1}|,|B_{2}|)$.
If for some $i\in \{ 1,2\}$, 
$|{ B_{i} \choose n} /E_{n}|=1$, then we put 
that all $n$-tuples $\bar{c}\in B_{1}\cup B_{2}$
meeting $B_{i}$ are pairwise $E_{n}$-equivalent.
We additionally arrange that they are equivalent
to some tuple from $B_{3-i}$ if $n\le |B_{3-i}|$. 
If $n\ge max(|B_{1}|,|B_{2}|)$ then all 
$n$-element $n$-tuples 
from $C$ are pairwise $E_{n}$-equivalent.
We take $E_n$ to be the minimal equivalence 
relation satisfying the conditions above. 
In particular if $n$-tuples $\bar{b}_1$ and $\bar{b}_2$ 
are $E_n$-equivalent to the same $n$-tuple from $A$, 
then $E_n (\bar{b}_1 ,\bar{b}_2 )$. 
\parskip0pt

We can now define the linear orderings $<_n$ 
on $C/E_n$ for $n\in P$.
There is nothing to do if $|{C \choose n} /E_{n}|=1$.
In the case when for some $i=1,2$,
$|{B_{i} \choose n} /E_{n}| = 1$, the relation $<_n$ 
is defined by its restriction to $B_{3-i}$.
When 
$|{B_{1} \choose n} /E_{n}|\not= 1\not= |{B_{2} \choose n} /E_{n}|$ 
and $V_1$, $V_2 $ is a pair of two 
$<_n$-neighbours among $E_n$-classes 
having representatives both in
${B_1 \choose n}$ and ${B_2 \choose n}$, 
we amalgamate the $<_n$-linear orderings 
between $V_1$ and $V_2$ assuming 
that all elements of 
${B_1 \choose n}/E_n \cap [V_1 ,V_2 ]$ 
are less than those from 
${B_2 \choose n}/E_n \cap [V_1 ,V_2 ]$. 

We appropriately modify this procedure 
for intervals open from one side. 
It is clear that this defines $<_n$-ordering 
on ${C \choose n} /E_{n}$. 
\parskip0pt

(2) The statement that $Th(M)$ admits elimination 
of quantifiers and is $\omega$-categorical, follows 
from (1). 
This also implies that $M$ is a natural 
expansion of $M_0$. 

To see the second statement of this 
part of the lemma 
%%%%%%%%%%%%%%%%%%%
%we apply the proof of Lemma \ref{support}. 
it is enough to show that for 
$n\in P$ and any two sequences 
$V_1 <_n ...<_n V_k$ and 
$V'_1 <_n ...<_n V'_k$ from $M_n$  
there is an automorphism of $M$ taking 
each $V_i$ to $V'_{i}$. 
To see this we use the fact that $M$ is 
the Fra\"{i}ss\'{e} limit of $\mathcal{K}_P$. 
This allows us to find pairwise disjoint 
representatives of classes $V_1 ,...,V_k$, 
say $\bar{a}_1 ,...,\bar{a}_k$,  
and classes $V'_1 ,...,V'_k$, 
say $\bar{a}'_1 ,...,\bar{a}'_k$, so that 
for every $m\not= n$ all $m$-tuples of the   
substructures $\bar{a}_1 \cup ...\cup \bar{a}_k$
and $\bar{a}'_1 \cup ...\cup \bar{a}'_k$ 
are $E_m$-equivalent. 
Moreover all $n$-tuples meeting at least 
two $\bar{a}_s$, $\bar{a}_t$ 
or $\bar{a}'_s$, $\bar{a}'_t$
also belong to a single $E_n$-class. 
Taking an appropriate isomorphism 
induced by these representatives 
we extend it to a required automorphism.  

%%%%%%%%%%%%%%%%%%%%%
%We now prove a helpful lemma. 
%
%
%\begin{lem} \label{support} 
%Let $\rho_i \in Sym (M_0 /E_i )$, $i\le k$, 
%be a sequence of finitary permutations. 
%Then there is an automorphism $\alpha \in Aut(M_0 )$ 
%realising each $\rho_i$ on its support.   
%\end{lem} 
%
%{\em Proof.} 
%%%%%%%%%%%%%%%%%% 
(3) We develop the argument of (2). 
For each $\rho_i$ find a sequence 
$\bar{a}_1 ,...,\bar{a}_t$ of pairwise disjoint 
tuples from $M$ representing the $E_i$-classes 
of the domain and of the range of $\rho_i$. 
We may assume that for any $j\not= i$ 
all $j$-tuples of  the union 
$\Omega_i = \bar{a}_1 \cup ...\cup \bar{a}_t$
belong to the same $E_j$-class. 
Moreover all $i$-tuples meeting at least two 
$\bar{a}_l$, $\bar{a}_m$  also form a single 
$E_i$-class. 
Thus $\rho_i$ can be realised by 
a  partial map on $\Omega_i$. 
We may arrange that all $\Omega_i$ are pairwise 
disjont and do not have common $E_n$-classes. 
Thus all $\rho_i$ can be realised by a partial 
isomorphism on the union of these $\Omega_i$. 
Since $M$ is ultrahomogeneous, this partial 
isomorphism can be extended to an 
automorphism of $M$. 
$\Box$ 
%%%%%%%%%%%%%%%%%%%
%\bigskip 
%
%It is clear that the statement of this lemma holds 
%if each $\rho_i$ is considered as a finite 
%permutation of some finite set (of classes) 
%containing its support instead the support in the formulation. 
%%%%%%%%%%%%%%%%%%
%$\Box$ 

\bigskip 

Let us consider $M$ in the language $L^S_P$, i.e. 
$$ 
(M , E_1 ,...,E_n ,...) \cup (M_1 , *_1 )\cup ...\cup (M_n ,*_n ) \cup ...\mbox{ , } 
$$ 
where $*_n = <_n$ for $n\in P$ and disappears 
for $n\not\in P$. 
By Lemma \ref{LEM}(3) the structure of all sorts  
$\{ M_n :n\in \omega\}$ coincides with the theory 
$\tilde{T}_1$ of sorts $\{ S_n :n\in \omega\}$ 
of the theory $T^S_P$. 
This implies the following corollary.

\begin{cor} \label{COR} 
The theory of $M$ in the language $L^S_P$ 
coincides with $T^S_P$. 
In particular the  sets $M_n$ are stably 
embedded into $M$. 
\end{cor} 

%\bigskip 

We see that for $n\in P$ 
any automorphism of $(M_n ,<_n )$ 
can be realized by 
an automorphism of $M$. 
Assume that 
$2n\not\in P$.  
Let us consider automorphisms 
$\alpha$ of $M_n$ which are 
{\em increasing}, i.e. for any 
$V\in M_n$, $V<_n \alpha(V)$. 

Take an orbit of $\alpha$  
of the following form:    
$$
... \rightarrow \bar{a}_{-1} \rightarrow \bar{a}_0 \rightarrow 
\bar{a}_1 \rightarrow \bar{a}_2 \rightarrow \bar{a}_3 \rightarrow \bar{a}_4 \rightarrow ...   
$$ 
and consider $E_{2n}$-classes of 
tuples $\bar{a}_i \bar{a}_{i+1}$.  
Applying ultrahomogenity and the choice 
of $n$ it is easy to see that $\alpha$ can be taken 
so that there are four $E_{2n}$-classes, say 
$V_1, V_2 , V_3 ,V_4$, represented by 
consecutive pairs of tuples 
$\bar{a}_1 ,\bar{a}_2 ,\bar{a}_3 ,\bar{a}_4 ,\bar{a}_5 ,\bar{a}_6$ 
and  $\alpha$ acts 
on them by $\mathbb{Z}/4\mathbb{Z}$: 
$$ 
\mbox{ if } \bar{a}_1 \bar{a}_2 \in V_1 
\mbox{ , then } 
\bar{a}_2 \bar{a}_3 \in V_2 \mbox{ , } 
\bar{a}_3 \bar{a}_4 \in V_3 
\mbox{ and } \bar{a}_4 \bar{a}_5 \in V_4 ,
$$  
where $\bar{a}_1 \bar{a}_2$   
and $\bar{a}_5\bar{a}_6$  
are $E_{2n}$-equivalent.  

Slightly generalising this situation 
we will say that a sequence 
$\bar{a}_1 ,\bar{a}_2 ,\bar{a}_3 ,\bar{a}_4 ,\bar{a}_5 ,\bar{a}_6$ 
is $<_n$-{\bf increasing of type} 
$\mathbb{Z}/4\mathbb{Z}$ 
if the following conditions are satisfied: 
\begin{itemize} 
\item tuples  $\bar{a}_1 \bar{a}_2$,  $\bar{a}_2 \bar{a}_3 $ and $\bar{a}_3 \bar{a}_4$ 
are of the same isomorphism type,  
\item tuples $\bar{a}_1 \bar{a}_2 \bar{a}_3 \bar{a}_4$ and $\bar{a}_3 \bar{a}_4 \bar{a}_5 \bar{a}_6$
are of the same isomorphism type and 
\item tupes $\bar{a}_1 \bar{a}_2$ and $\bar{a}_5 \bar{a}_6$ are 
$E_{2n}$-equivalent but not $E_{2n}$-equivalent to $\bar{a}_3 \bar{a}_4$. 
\end{itemize}

Let $L'$ be an extension of $L_P$ and 
$M' =(M,\bar{{\bf r}})$ be an $L'$-expansion 
of $M$ with quantifier elimination. 
We do not demand that $\bar{{\bf r}}$ 
is finite, we only assume that $M'$ is 
a precompact expansion.  
It is clear that $M'$ induces 
a subgroup of $Aut(M_n ,<_n )$.

We will say that a sequence 
$\bar{a}_1 ,\bar{a}_2 ,\bar{a}_3 ,\bar{a}_4 ,\bar{a}_5 ,\bar{a}_6$ 
is $<_n$-{\bf increasing of type} 
$\mathbb{Z}/4\mathbb{Z}$ {\bf  in} $M'$ 
if the definition above holds under the assumption 
that the isomorphism types appeared in 
the definition are considered with 
respect to the relations of $M'$. 

\begin{thm} \label{THM} 
Let $M$ be the generic structure of 
${\cal K}_P$ where $P\not=\emptyset$. 
Then the group $G=Aut(M)$ is amenable, 
$M$ does not satisfy Hrushovski's extension 
property and does not have an extremely 
amenable ultrahomogeneous expansion 
by a linear ordering. 

Let $M'$ be a precompact expansion of 
$M$ with quantifier elimination. 
If $Aut(M')$ is extremely amenable, 
then for any $n\in P$ with $2n \not\in P$ 
the structure $M'$ does not have 
an $<_n$-increasing sequence  
of type $\mathbb{Z}/4\mathbb{Z}$. 
\end{thm} 

The main point of this theorem is 
that although in different arities 
the structures induced by $M$ 
are completely independent, 
any expansion $M'$ as in the formulation 
simultaneously destroys $M$ 
in all arities $n\in P$ with $2n\not\in P$. 

The proof below uses the proof of 
Theorem \ref{nR}. 
\bigskip 

{\em Proof of Theorem \ref{THM}.} 
For each $n>1$ enumerate all $E_n$-classes. 
Consider the expansion of $M$ by distinguishing 
each $E_n$-class by a predicate $P_{n,i}$ according 
the enumeration. 
Let $L^*$ be the language of all predicates $P_{n,i}$ 
and let $M^*$ be the $L^*$-structure defined on $M$. 
By Claims 1 - 4 of the proof of Theorem \ref{nR} 
the structure $M^*$ has Hrushovski's 
extension property and $Aut(M^* )$ 
is amenable. 

Let us consider the structure 
$(M_n , <_n )$, where $n\in P$. 
As it is isomorphic to $(\mathbb{Q}, <)$, 
the group $Aut (M_n , <_n )$ 
is extremely amenable (\cite{KPT}). 

Since each automorphism of $M$ 
preserves all $<_i$, $i\in P$, 
it is easy to see that there is a natural 
homomorphism from $Aut(M)$ to the product 
$$
\prod_{i\in P} Aut (M_i , <_i ) \times   
\prod_{i\not\in P} Sym (M_i )
$$  
and $Aut(M^* )$ is the kernel of it. 
By Corollary \ref{COR}  
this homomorphism is surjective. 
Now by Theorem 449C of \cite{fremlin} 
we have the following claim. 

{\em The group $Aut(M)$ is amenable.}  

To see that $M$ does not satisfy 
Hrushovski's extension property 
take $n\in P$ and let us 
consider any triple of pairwise disjoint 
$n$-tuples $\bar{a}$, $\bar{b}$, $\bar{c}$ 
representing pairwise 
distinct elements of $M_n$ so that 
$$ 
\bar{a} <_n \bar{b} <_n \bar{c}. 
$$ 
Then the map $\phi$ fixing $\bar{a}$ and 
taking $\bar{b}$ to $\bar{c}$ cannot 
be extended to an automorphism 
of a finite substructure of $M$. 

Consider a linearly ordered expansion 
$(M,<)$ with quantifier elimination.  
To see that $Aut(M,<)$ is 
not extremely amenable just apply the 
argument of statement (b) of Theorem \ref{nR}.  
Since at arity 2 the structure $M$ 
coincides with $M_0$ it works 
without any change. 

To prove the second part of the theorem 
we slightly modify that argument. 

Let $n\in P$ and  
$2n \not\in  P$. 
Let a structure $B$ consist of $6n$ 
elements forming a sequence   
$$
\bar{a}_1 <_n \bar{a}_2 <_n \bar{a}_3 <_n \bar{a}_4 <_n \bar{b}_1 <_n \bar{b}_2  , 
$$ 
where the tuples $\bar{a}_1 \bar{a}_2 $ and 
$\bar{b}_1 \bar{b}_2$ are $E_{2n}$-equivalent 
but not of the same $E_{2n}$-class 
with $\bar{a}_3 \bar{a}_4$.  
We assume that the tuples  
$\bar{a}_1 \bar{a}_2$, 
$\bar{a}_2 \bar{a}_3$, 
and $\bar{a}_3 \bar{a}_4$
are of the same isomorphism class in $M'$ 
and  the substructure  
$\bar{a}_1 \bar{a}_2 \bar{a}_3 \bar{a}_4 <M'$ 
is isomorphic to 
$\bar{a}_3 \bar{a}_4\bar{b}_1 \bar{b}_2<M'$. 
Since $Aut(M')$ is extremely amenable, these structures 
are rigid and the corresponding isomorphisms are uniquely 
defined on these tuples. 

Let $A$ represent the isomorphism 
class of 
$\bar{a}_1 \bar{a}_2 \bar{a}_3 \bar{a}_4$ 
in $M'$.  
Let us show that the Ramsey property 
does not hold for the age of $M'$.  
Take any finite substructure $C$ of this 
age which extends $B$. 
Fix any enumeration of $E_{2n}$-classes 
ocurring in $C$. 
Then colour a copy of $A$ red if 
the class of the first two $n$-tuples  
is enumerated before the class 
of the last pair. 
Otherwise colour such a copy green. 
It is clear that $C$ does not contain 
a structure isomorphic to $B$ so that all 
substructures of  type $A$ are of the same colour. 
$\Box$

%%%%%%%%%%%%%%%%%%%%%%%%%%%%%%%%%%%

%\bigskip 

%\newpage

\bigskip

Institute of Mathematics, University of Wroc{\l}aw, pl.Grunwaldzki 2/4, 50-384 Wroc{\l}aw, Poland, \\
 E-mail: ivanov@math.uni.wroc.pl 

\end{document}